\documentclass[preprint,3p,11pt,sort&compress]{elsarticle}

 \usepackage{hyperref}
 \hypersetup{backref, colorlinks=true, linkcolor=blue}

\usepackage{stfloats}
\usepackage{subcaption}

 \usepackage{pict2e}

\usepackage{graphicx} 
\usepackage{booktabs} 
\usepackage{microtype} 
\usepackage{amssymb,amsmath}
\usepackage[ruled,vlined]{algorithm2e}

\begin{document}
\begin{center}
{\bf \Large Multilevel Method 
with Low-Order Equations of Mixed Types 
and  \\ \vspace{0.1cm}
Two Grids in Photon Energy
  for Thermal Radiative Transfer}
\end{center}

\author[ncsu]{Dmitriy Y. Anistratov}
\ead{anistratov@ncsu.edu}
\author[llnl]{Terry S. Haut}
\ead{haut3@llnl.gov}
\address[ncsu]{Department of Nuclear Engineering,
North Carolina State University Raleigh, NC}
\address[llnl]{Lawrence Livermore National Laboratory, Livermore, CA}

\begin{frontmatter}
\begin{abstract}
Thermal radiative transfer (TRT)   is an essential piece of physics  in  
 inertial confinement fusion, high-energy density physics,  astrophysics etc.
The physical models of this type of problem are defined  
  by strongly coupled differential  equations describing multiphysics phenomena.
This paper presents a new nonlinear  multilevel iterative method  with two photon energy grids
for solving the multigroup radiative transfer equation (RTE) coupled with the material energy balance equation (MEB).
The multilevel  system of equations of the method is formulated by means of a nonlinear projection approach.
The  RTE is projected over elements of phase space to derive the low-order equations of different types.
The hierarchy of equations consists of
(1) multigroup weighted flux equations which can be interpreted as the multigroup RTE averaged over
subintervals of angular range and (2) the effective grey (one-group) equations which are spectrum averaged 
low-order quasidiffusion (aka variable Eddington factor) equations.
The system of  RTE, low-order and MEB equations is approximated by the
fully implicit Euler time-integration method in which absorption coefficient   and emission term are evaluated
at the current time step.
Numerical results are presented to demonstrate   convergence of a multilevel iteration algorithm  
 in the Fleck-Cummings test problem with Marshak wave solved with
  large number of photon energy groups. 
\end{abstract}
\begin{keyword}
high-energy-density physics,
radiative transfer,
multilevel iteration methods,
multigrid methods,
multiphysics problems.
\end{keyword}
\end{frontmatter}

\section{Introduction}\label{sec:1}

Radiative transfer plays an important part in   high-temperature phenomena of
high-energy density physics, inertial confinement fusion, and astrophysics \cite{drake-2006}.
High-energy photons   affect  redistribution of energy in such physical systems.
The state of the physical system  is significantly impacted  by absorption of photons by matter and their emission, 
leading to change in  the material energy.  
The Boltzmann  radiative transfer   equation (RTE) describes photon transport in a medium and  the interaction of particles  with matter.
The absorption coefficient and emission source term in the RTE nonlinearly depend on  the material temperature.

The fundamental thermal radiative transfer (TRT) model in 1D slab geometry is defined by
the time-dependent multigroup  RTE for
the group intensity
\begin{equation} \label{rt}
\frac{1}{c}  \partial_t I_g  + \mathcal{L} I_g = \sigma_{g} B_g \, , \quad
g \in \mathbb{N}(n_{\nu}) \, ,
\end{equation}
\begin{equation}
 \mathcal{L} I_g(x,\mu,t) \equiv \mu  \partial_x I_g(x,\mu,t)   + \sigma_{g}(x;T)I_g(x,\mu,t) \, ,
\end{equation}
\begin{equation}
I_g\big|_{t=0} = I_g^0 \, , \quad I_g\Big|_{ \overset{\mu>0}{x=0} }  = I_g^+  \, , \quad I_g\Big|_{ \overset{\mu<0}{x=X} }  = I_g^- \, ,
\end{equation}
and the material energy balance (MEB) equation
\begin{equation}\label{eb}
\ \partial_t \varepsilon(T)  =
 \sum_{g=1}^{n_{\nu}}\int_{-1}^1 \sigma_{g}(x;T)\Bigl( I_{g}(x,\mu,t) -  B_{g}(T)\Bigr) d \mu   \, ,
\end{equation}
\begin{equation}
T\big|_{t=0} = T^0 \, .
\end{equation}
Here $x$ is the spatial position;
$\mu$ is     the directional cosine of particle motion;
t is time;
$g$ is the index of the photon energy group;
 $T$ is the material temperature;
$\varepsilon$ is the material energy density;
\begin{equation}
I_g = \int_{\nu_g}^{\nu_{g+1}} I_{\nu} d \nu \,  
\end{equation} 
is the group  intensity,  $\nu$ is the photon frequency,
\begin{equation}
B_g(T)=\int_{\nu_g}^{\nu_{g+1}} B_{\nu}(T) d \nu \,  
\end{equation} 
is the group Planckian function, 
\begin{equation}
\sigma_g(T)=\frac{\int_{\nu_g}^{\nu_{g+1}} \sigma_{\nu}(T) B_{\nu}(T) d \nu}{ \int_{\nu_g}^{\nu_{g+1}}   B_{\nu}(T) d \nu} \,  
\end{equation}
is the group averaged  absorption coefficient.
The group structure is defined by the photon energy grid
\begin{equation}\label{fine-grid}
\Lambda_{h\nu}=\{h\nu_g, g \in \mathbb{N}(n_{\nu}+1): \ h\nu_{min} =  h\nu_1 <   \ldots  <  h\nu_g   <   h\nu_{g+1} <  \ldots  <   h\nu_{n_{\nu}+1}  =   h\nu_{max}\} \, .
\end{equation}
 This  TRT model   neglects material motion, scattering, heat conduction, and external sources.
 It describes the phenomena of a supersonic radiation wave.
The general model is formulated by the system of radiation hydrodynamics equations, which involves conservation of mass and momentum equations and accounts for the radiation effect on matter \cite{zel-1966}.

There exist a family of computational methods for TRT problems  which  apply 
linearizaton to the  multigroup RTE resulting 
in the RTE with a  pseudo-scattering  term.  
To accelerate transport iterations due to   pseudo-scattering,
efficient  algorithms are used
\cite{morel-jqsrt-1985,larsen-jcp-1988,morel-jcp-2007,mla-m&c2013}.
A different group of methods follows the moment-based approach to formulate   iteration schemes  based  high-order and low-order equations \cite{park-et-al-2012,yee-2016,abw-2017}.

In this paper, we apply a nonlinear projection approach  to formulate a new multilevel iteration method for solving coupled RTE and MEB equation  \cite{gol'din-cmmp-1964,Goldin-sbornik-82,dya-vyag-ttsp,dya-jcp-2021,sam-et-al-jcp-2023,sam-terry-jctt-2023}.
The hierarchy of equations of the method is defined for  different moments of the radiation intensity.
The  RTE is projected sequentially over elements of phase space to derive  nonlinear  moment equations with exact closures.
One part of the hierarchy of moment equations  is defined by multigroup low-order nonlinear weighted flux (LONWF) equations,
 which can be viewed as the  RTE averaged over subintervals of angular range with certain weight functions
 \cite{dya-vyag-ttsp,germ-DAN-68,goldin-69,surzh-jep-82,dya-ewl-jcp-2001,lr-dya-ttsp-2007,lr-dya-2010,shands-hanopy-morel-m&c2023}.
Another part of the hierarchy is the effective grey  
low-order quasidiffusion (GLOQD), aka variable Eddington factor (VEF), equations  \cite{gol'din-cmmp-1964,auer-mihalas-1970,PASE-1986,Winkler-85,dya-aristova-vya-mm1996,aristova-vya-avk-m&c1999}.
The LONWF equations are defined for partial angular  moments of the intensity. 
They mimic the structure of the RTE and lead to a coupled set of Cauchy problems in space.
The GLOQD equations are formulated for the full-range angular moments of the intensity  and have a structure similar to the $P_1$ equations,
resulting in  a boundary value problem in space. Thus, the proposed multilevel method is based on low-order problems of mixed types.
Parallel algorithms developed for solving the RTE on GPUs can be adopted for solving the LONWF equations.
This   provides an   opportunity for developing fast and efficient algorithms for TRT problems on advanced
high-performance computing system with GPUs.

Our approach possesses the following major potential advantages over more
traditional discretization schemes:
(i) the low-order system for the radiation energy is hyperbolic and this naturally avoids some pathological issues
 that can arrive when coupling to diffusion near highly heterogeneous material interfaces, 
(ii) our formulation can result in positivity for both the angular intensity and the angle-averaged solutions through positivity-preserving sweeps, (iii) 
our method naturally extends to using a hierarchy of levels in both angle and frequency, which opens the way for potentially much more efficient iterative methods, 
(iv) our approach can dynamically adapt the multi-level hierarchy of equations for different physics regimes, potentially resulting in greater efficiency.

The remainder of this paper is organized as follows.
In Sec. \ref{method},  the multilevel method is formulated.
The iteration algorithm  is  described in Sec.  \ref{algorithm}.
In  Sec. \ref{results}, the numerical results are presented.
We conclude with a brief summary in Sec. \ref{conclusion}.

\section{Formulation of the Multilevel Method \label{method}}

The   method is formulated by a multiscale system of differential equations consisting of
\begin{itemize}
\item  the multigroup high-order RTE  (Eq. \eqref{rt}) on the grid $\Lambda_{h\nu}$ for the group intensities $I_g$, 
\item  the multigroup LONWF equations on the grid $\Lambda_{h\nu}$
for the  group partial angular moments of radiation   intensities
\begin{equation}
\phi_g^{\pm}  =\pm  \int_0^{\pm 1}  I_g   d \mu \, , 
\end{equation}
\item  effective    GLOQD  equations on the coarse photon energy grid
 \begin{equation} 
\Lambda_{h\nu}^{\ast} =\{ h\nu_1 = h\nu_{min}, h\nu_2 = h\nu_{max} \} \, 
 \end{equation}
 for the total radiation energy density  and flux 
 \begin{equation} \label{E&F}
E= \sum_{g=1}^{n_{\nu}} E_g \,  , \quad  F= \sum_{g=1}^{n_{\nu}} F_g \,  ,
\end{equation}
where the group quantities are defined by
 \begin{equation} \label{Eg&Fg}
 E_g = \frac{1}{c}  \int_{-1}^{1}  I_g   d \mu \,   , \quad  F_g =   \int_{-1}^{1} \mu  I_g   d \mu \, .
\end{equation}

\end{itemize}
The multigroup  LONWF equations for $\phi_g^{\pm}$  are derived by taking angular moments
of the RTE with weight $w(\mu) = \gamma (1 + \sqrt{3}|\mu|)$ over half angular ranges, and are given by  \cite{dya-ewl-jcp-2001,lr-dya-ttsp-2007}
\begin{equation} \label{f-lo-eqs}
 \frac{1}{c}\partial_t (K_g^{\pm} \phi_g^{\pm})
\pm  \partial_x (H_g^{\pm} \phi_g^{\pm})  
+  \sigma_{g}  K_g^{\pm}\phi_g^{\pm}
=   \sigma_{g}B_g \, ,
\end{equation}
\begin{equation}
\gamma = \Big( 1 + \frac{1}{2}\sqrt{3} \Big)^{-1} \, .
\end{equation}
The closures for the low-order equations \eqref{f-lo-eqs} are formulated by the factors
\begin{equation}
K_g^{\pm} = \gamma \frac{\int_0^{\pm 1} \big(1 + \sqrt{3} |\mu|I_g \big) d \mu  }{\int_0^{\pm 1}I_g  d \mu} \, ,
\quad 
H_g^{\pm} = \gamma \frac{\int_0^{\pm 1}|\mu| \big(1 + \sqrt{3} |\mu|I_g \big) d \mu  }{\int_0^{\pm 1}I_g  d \mu} \, .
\end{equation}
The group radiation energy density and flux (Eq. \eqref{Eg&Fg}) are defined with the LONWF solution as follows:
\begin{equation}
E_g = \frac{1}{c}(\phi_g^- +   \phi_g^+) \, , \quad
F_g = G_g^-\phi_g^- + G_g^+\phi_g^+  \, ,  
\end{equation}
where
\begin{equation}
G_g^{\pm}  = \frac{\int_0^{\pm 1}  \mu  I_g d \mu} {\int_0^{\pm 1}   I_g d \mu } \, .
\end{equation}
The projection with the angular weight  $w(\mu)$ leads to the low-order equations which meet the asymptotic diffusion limit in the interior of a diffusive domain with asymptotic boundary conditions.
The low-order equations \eqref{f-lo-eqs} and high-order RTE can be approximated by independent discretization schemes.

The first step in deriving the effective GLOQD  equations is integrating the RTE over $-1 \le \mu \le 1$ with
weigh 1 and $\mu$ and introducing exact closure by means of the group Eddington factor
\begin{equation}
 f_g   = \frac{\int_{-1}^{1}  \mu^2  I_g d \mu} {\int_{-1}^{1}   I_g d \mu }  \, 
\end{equation}
to obtain multigroup moment equations for $E_g$ and $F_g$. 
On the next step, the moment equations are summed over photon energy groups to get  
the effective grey moment  equations  on the coarse photon energy grid $\Lambda_{h\nu}^{\ast}$ with just one group ($n_{\nu}=1$)
 for the  total radiation energy density   and flux  (Eq. \eqref{E&F}).
The effective GLOQD equations  are defined by
 \cite{PASE-1986,dya-aristova-vya-mm1996,aristova-vya-avk-m&c1999,dya-jcp-2019}:
\begin{subequations}\label{gloqd}
\begin{equation}\label{gloqd-1}
 \partial_t E 
+  \partial_x F 
+ c \bar \sigma_E  E
= c  \bar \sigma_B a_R T^4 \, ,
\end{equation}
\begin{equation}\label{gloqd-2}
\frac{1}{c} \partial_t F 
+ c \partial_x (\bar f E)
+ \bar \sigma_R  F + \bar \eta E = 0 \, ,
\end{equation}
\end{subequations}
where
\begin{equation}
\bar \sigma_E= \frac{\sum_{g=1}^{n_{\nu}} \sigma_{g}E_g}
{\sum_{g=1}^{n_{\nu}} E_g}  \, ,
\quad
\bar \sigma_{B}= \frac{\sum_{g=1}^{n_{\nu}}\sigma_{g} B_g}
{\sum_{g=1}^{n_{\nu}} B_g}  \, ,
\quad
\bar \sigma_R= \frac{\sum_{g=1}^{n_{\nu}} \sigma_{g}|F_g|}
{\sum_{g=1}^{G} |F_g|}  \,
\end{equation}
are different grey absorption coefficients.
\begin{equation}
\bar f = \frac{\sum_{g=1}^{n_{\nu}} f_g E_g}
{\sum_{g=1}^{n_{\nu}} E_g}  \,  
\end{equation}
is  the grey Eddington factor. It is a spectrum average of group Eddington factors.
The  compensation term in Eq. \eqref{gloqd-2} is defined by
\begin{equation}
\bar \eta  = \frac{\sum_{g=1}^{n_{\nu}} (\sigma_{g} - \bar \sigma_R) F_g}
{\sum_{g=1}^{n_{\nu}} E_g} \, .
\end{equation}
The MEB equation in the effective grey form  is given by
\begin{equation}\label{eb-2}
 \partial_t \varepsilon(T)  =  c  \big(\bar \sigma_E E - \bar\sigma_B a_R T^4\big) \, .
\end{equation}

Finally, the multilevel method  is defined by the hierarchy of  equations
 the general form of which is the following:
\begin{itemize}
\item the multigroup RTE
\begin{equation} \label{rt-geb}
\frac{1}{c}  \partial_t I_g  + \mathcal{L} I_g = q_g \, , \quad
g \in \mathbb{N}(n_{\nu}) \, ,
\end{equation}
\begin{equation*} 
\mathcal{L} = \mathcal{L}(T) \, , \quad  q_g =q_g (T) \, ,
\end{equation*} 
\item the  multigroup LONWF equations
\begin{equation} \label{group-lonwf-gen}
 \partial_t  (\mathcal{K}\boldsymbol{\Phi}_g)  +  \mathcal{P} \boldsymbol{\Phi}_g = \mathbf{Q}_g \, ,
\quad
\boldsymbol{\Phi}_g = ( \phi_g^-,  \phi_g^+ )^{\top} \, , \quad \mathbf{Q}_g = \mathbf{Q}_g(T) \, ,
\end{equation}
\begin{equation*} 
 \mathcal{P} =  \mathcal{P}(T) \, , 
\end{equation*} 
\item the GLOQD equations
\begin{equation} \label{grey-loqd-gen}
 \partial_t \mathbf{Y}  +  \mathcal{M} \mathbf{Y} = \mathbf{Q} \, ,
\quad
\mathbf{Y} = ( E,  F )^{\top} \, ,
\end{equation}
\begin{equation*} 
\mathcal{M} = \mathcal{M}(\boldsymbol{\Phi},T) \, , \quad  \mathbf{Q} =\mathbf{Q} (T) \, , \quad
\boldsymbol{\Phi} = \{\boldsymbol{\Phi}_g\}_{g=1}^{n_{\nu}}
\end{equation*}  
\item the MEB equation
\begin{equation}\label{grey-meb-gen}
 \partial_t \varepsilon  =  \mathcal{R}  \mathbf{Y}  - Q_{pl} \, ,
\end{equation}
\begin{equation*} 
 \mathcal{R} = \mathcal{R}(\boldsymbol{\Phi},T) \, , \quad Q_{pl} =  Q_{pl}(T) \, .
\end{equation*} 
\end{itemize}

In this study, we discretize this system of equations with a fully-implicit temporal scheme based on the backward Euler time-integration method.
The RTE is approximated in space by the simple corner balance scheme  \cite{mla-ttsp-1997}. 
The lumped linear discontinuous (LLD) scheme is used  to discretize the LONWF equations \cite{ewl-jem-jcp-1989,dya-ewl-jcp-2001,lr-dya-ttsp-2007}.
The GLOQD equations are approximated by the second-order finite volume method \cite{dya-jcp-2019}.
We note that independent spatial discretization  schemes are applied for approximating the RTE, LONWF, and GLOQD equations \cite{dya-vyag-ttsp}.

\section{Iteration Algorithm  \label{algorithm}}

\begin{figure}[h!]
\centering
{
\setlength{\unitlength}{0.8cm}
\begin{picture}(5,5)
\put(0,5){\makebox(0,0)[c]{ \textrm{grid}} }
\put(0,4){\makebox(0,0)[c]{ $\Lambda_{h\nu}$} }
\put(0,0){\makebox(0,0)[c]{ $\Lambda_{h\nu}^{\ast}$} }
\put(2,4){\line(1.5,-4){1.5}}
\put(3.5,0){\line(1.5,4){1.5}}
\put(2,4){\circle*{0.3}}
\put(1.4,4){\makebox(0,0)[c]{ $\boldsymbol{\Phi}_g$} }
\put(3.5,0){\circle*{0.3}}
\put(3.0,0){\makebox(0,0)[c]{ \textsf{T}} }
\end{picture}
}
   \caption{Diagram of
  the \textsf{V}-cycle of inner iterations.
   $\boldsymbol{\Phi}_g$ - calculation  of spectrum functions on the fine  grid   by solving the multigroup LONWF equations,
 \textsf{T} - calculation of temperature by solving coupled GLOQD and MEB equations on the coarse grid.}
  \label{V-cycle}
\end{figure}
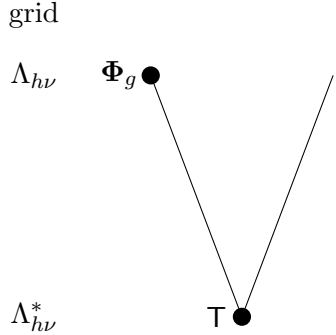

Algorithm \ref{V} presents  the multilevel method 
with    two photon energy grids for solving  the low-order equations \eqref{group-lonwf-gen} and \eqref{grey-loqd-gen} coupled with MEB equation \eqref{grey-meb-gen}.
On each time step, the outer iteration  is the transport iteration.
 $s$~is the   index of the outer iterations.
Note that there are no transport sweeps for $s=0$.
The inner iteration cycles consist  of  solving (i)   the system of the group LONWF equations 
on the given grid $\Lambda_{h\nu}$  and 
(ii)  the coupled GLOQD and MEB equations   on the coarse grid $\Lambda_{h\nu}^{\ast}$.
 $\ell$~is the   index of the inner iteration cycles.
 This inner iteration scheme is the multigrid $\textsf{V}$-cycle over two photon energy grids
\cite{dya-jcp-2021,multigrid}. Figure \ref{V-cycle} shows the diagram of the $\textsf{V}$-cycle.
The coupled  GLOQD and MEB  are linearized and use the Fr$\acute{\mbox{e}}$chet derivative   of   $\bar \sigma_E$
   to account for   variation  in group absorption coefficients and spectrum  due to change in temperature \cite{dya-jcp-2021}. 
On every cycle, the linearized system of the GLOQD and MEB equations is solved with one Newton iteration.

\medskip
\begin{algorithm}[H]
\DontPrintSemicolon
\While{$t^j \le t^{end}$}{
$s=0$, $T^{(0)}=T^{j-1}$, 
$K_g^{\pm \,  (1/2)}=K_g^{\pm \,  j-1}$\;
\While{$|| T^{(s)} - T^{(s-1)} ||>\epsilon||T^{(s)}|| \vee
 ||E^{(s)} - E^{(s-1)}||> \epsilon||E^{(s)}|| $ }{
\If{$s>0$}{
update group absorption coefficients $\sigma_g(T^{(s)})$\;
solve  $\boxed{c^{-1}  \partial_t I_g  + \mathcal{L} I_g = q_g}$ for $I_g^{(s+1/2)}$ on $\Lambda_{h\nu}$ for all $g$\;
compute group   factors  $K_g^{\pm \,  (s+1/2)}$, $H_g^{\pm \, (s+1/2)}$, $G_g^{\pm \, (s+1/2)}$, and $f_{g}^{(s+1/2)}$\;}
\While{$\ell \le \ell_{max} \wedge \Big( || T^{(\ell,s)} - T^{(\ell-1,s)} ||>\tilde\epsilon||T^{(\ell,s)}|| \vee
 || E^{(\ell,s)} - E^{(\ell-1,s)} ||>\tilde\epsilon||E^{(\ell,s)}|| \Big)$ }{
 given: $T^{(\ell,s)}$\;
update group absorption coefficients $\sigma_g(T^{(\ell,s)})$\;
solve   $\boxed{\partial _t (\mathcal{K}_g \boldsymbol{\Phi}_g) + \mathcal{P}_g \boldsymbol{\Phi}_g = \mathbf{Q}_g}$ 
for $\phi_g^{\pm \, (\ell,s)}$ 
on  $\Lambda_{h\nu}$  for all $g$\;
compute $E_g^{(\ell,s)}$    and $F_g^{(\ell,s)}$\; 
compute   $\bar \sigma_E^{(\ell,s)}$\!,  $\bar \sigma_{R}^{(\ell,s)}$\!,   $\bar{f}^{(\ell,s+1/2)}$\!,
$\bar{\eta}^{(\ell,s+1/2)}$\!,
and  Fr$\acute{\mbox{e}}$chet derivative $\delta_T \bar \sigma_E^{(\ell,s)}$ \;
solve  $\boxed{\partial _t \mathbf{Y}+ \mathcal{\bar M} \mathbf{Y} = \mathbf{Q} \quad  \partial_t \varepsilon   =  \mathcal{R} \mathbf{Y} - Q_{pl}}$ on $\Lambda_{h\nu}^{\ast}$   for $T^{(\ell+1,s)}$, $E^{(\ell+1,s)}$,   $F^{(\ell+1,s)}$\;
}
$T^{(s+1)} \leftarrow  T^{(\ell+1,s)}$, $s \leftarrow  s+1$ \;
}
}
\caption{\label{V} The multilevel method with two photon energy groups and $\textsf{V}$-cycle
}
\end{algorithm}

\section{Numerical Results \label{results}}

\begin{figure}[hbt]
\centering
\subfloat[Temperature   \label{T}]{\includegraphics[scale=0.3]{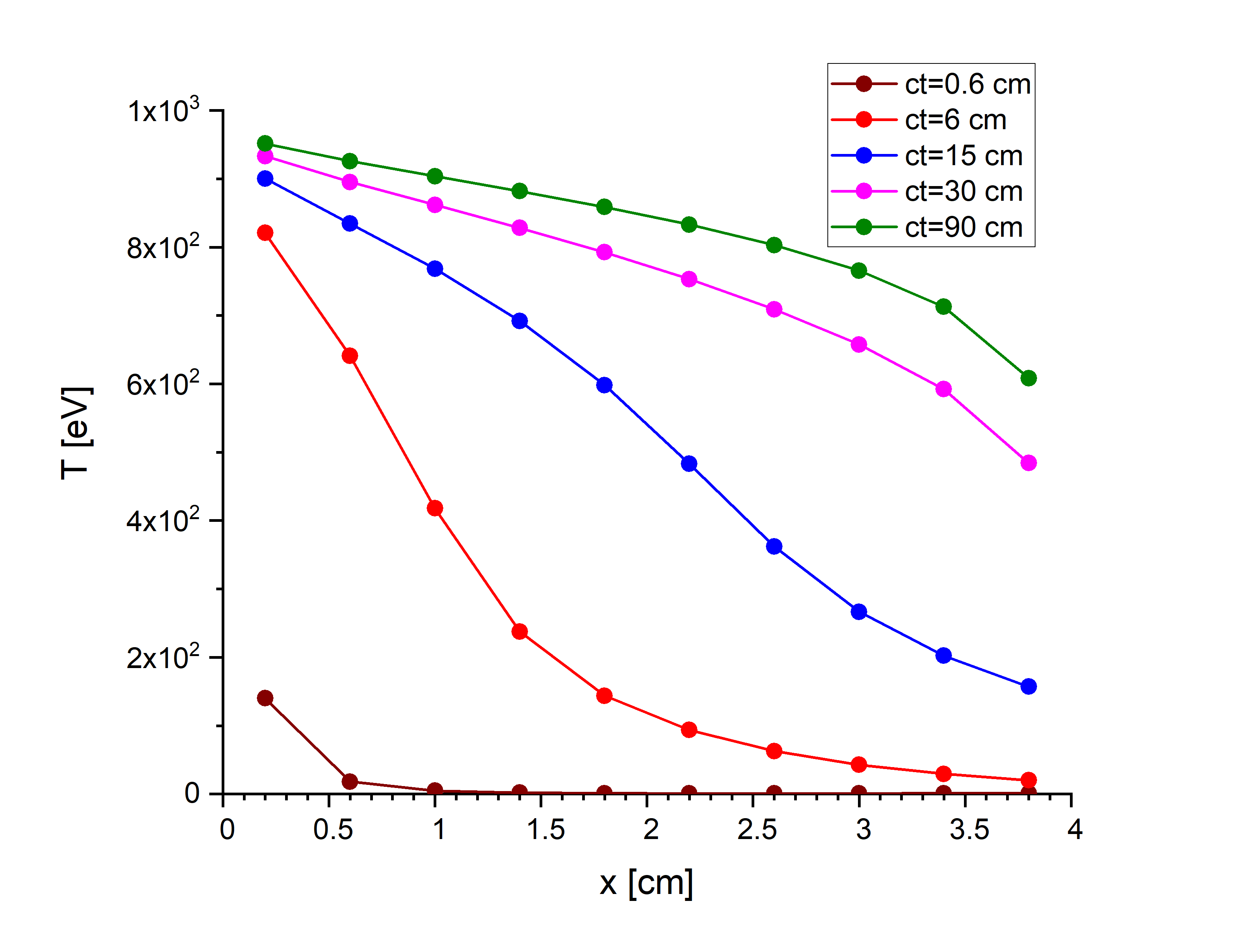}}
\subfloat[Energy density \label{E}]{\includegraphics[scale=0.3]{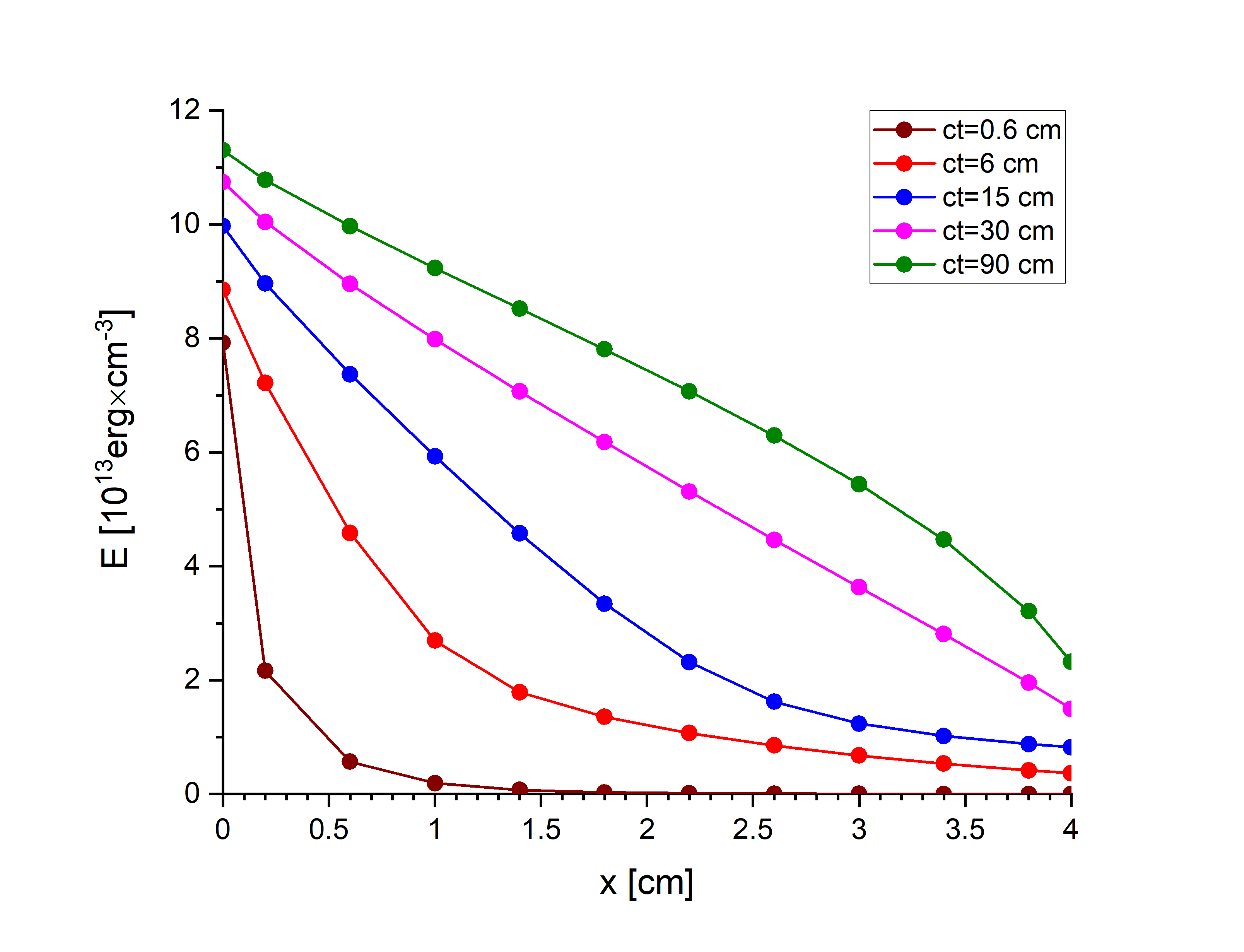}}
\caption{\label{solution}
Numerical solution of the FC test with 256 groups using $\Delta t=2 \! \times \! 10^{-3}$ sh.}
\end{figure}

The convergence of the presented iteration algorithm is analyzed on    the Fleck-Cummings (FC) test \cite{fleck-1971}.
It is defined for a  slab ($0\le x \le 4$ cm).
The  analytic spectral absorption coefficient of the material is given by
\begin{equation}
\sigma_{\nu}(T)    =    \frac{27}{(h\nu)^3}\big(1-e^{-\frac{h\nu}{kT}}\big) \, .
\end{equation}
The incoming radiation on the left boundary has the black-body spectrum at $kT_b$=1 ~keV and hence $I_{g}^+=B_g(T_b)$.
The right boundary is vacuum.
At $t=0$, the  temperature in the slab equals $kT_0$=10$^{-3}$ keV and  the group intensity at the initial instant is  $ I_g^0 =   B_{g}(T_0)$.
The material energy is a linear function and given by $\varepsilon(T)   =  c_v T$, where  $c_v  = \    0.5917a_R  T_b^3$.
The spatial  mesh is uniform with  10 intervals.
The time interval of the problem is  $t \in [0, 0.3 \, \mbox{sh}]$ ($1~sh~=~10^{-8}$~sec).
The photon energy grid consists of $n_{\nu}=256$   groups over the interval $ h\nu \in [0, 10^7]$~keV.
 $n_{\nu}-2$  groups are evenly spaced in logarithmic scale between $h\nu_a= 10^{-4}$~keV and $h\nu_b = 10$ keV.
The angular grid has  16  angular directions defined by the double S$_8$  Gauss-Legendre quadrature set.
The parameters of convergence criteria  are $\epsilon  = 10^{-6}$ and $\tilde \epsilon  = 10^{-7}$.
Figure~\ref{solution} shows temperature and total energy density at various instants of time obtained with $\Delta t = 2 \times 10^{-3}$ sh and  illustrates evolution of    heat and radiation waves.

Table   \ref{tbl-256-004} presents the total numbers of transport iterations, $N_{ti}$, and
   cycles, $N_c$,  accumulated over all time steps when the problems is solved with  
    $\Delta t=1 \! \times \! 10^{-3}, 2 \! \times \! 10^{-3}, 4 \! \times \! 10^{-3}$ sh.
    \begin{itemize}
    \item  $\mathbf{\Delta t=1 \! \times \! 10^{-3} \, sh}$. 
     The algorithm with 4 cycles per transport iteration 
    yields the minimum number of transport iterations ($N_{ti}=653$) over the time range of the problem.
    The algorithm with 3 cycles per transport iteration leads to extra 5 transport iterations with smaller
    total number of cycles. These two algorithms have comparable computational costs.
        \item  $\mathbf{\Delta t=2 \! \times \! 10^{-3} \, sh}$.   
        The algorithm with 4 cycles per transport iteration  has
             the optimum number of transport iterations ($N_{ti}=384$). The algorithm with 3 cycles needs extra 15 transport iterations. 
            \item  $\mathbf{\Delta t=4 \! \times \! 10^{-3} \, sh}$.     
                 The algorithm with 5 cycles per transport iteration  requires
             the minimum number of transport iterations ($N_{ti}=232$).
    \end{itemize}

\begin{figure}[h!]
\centering
\subfloat[$\ell_{max}$=3  \label{V3-1e-3}]
{\includegraphics[scale=0.28]{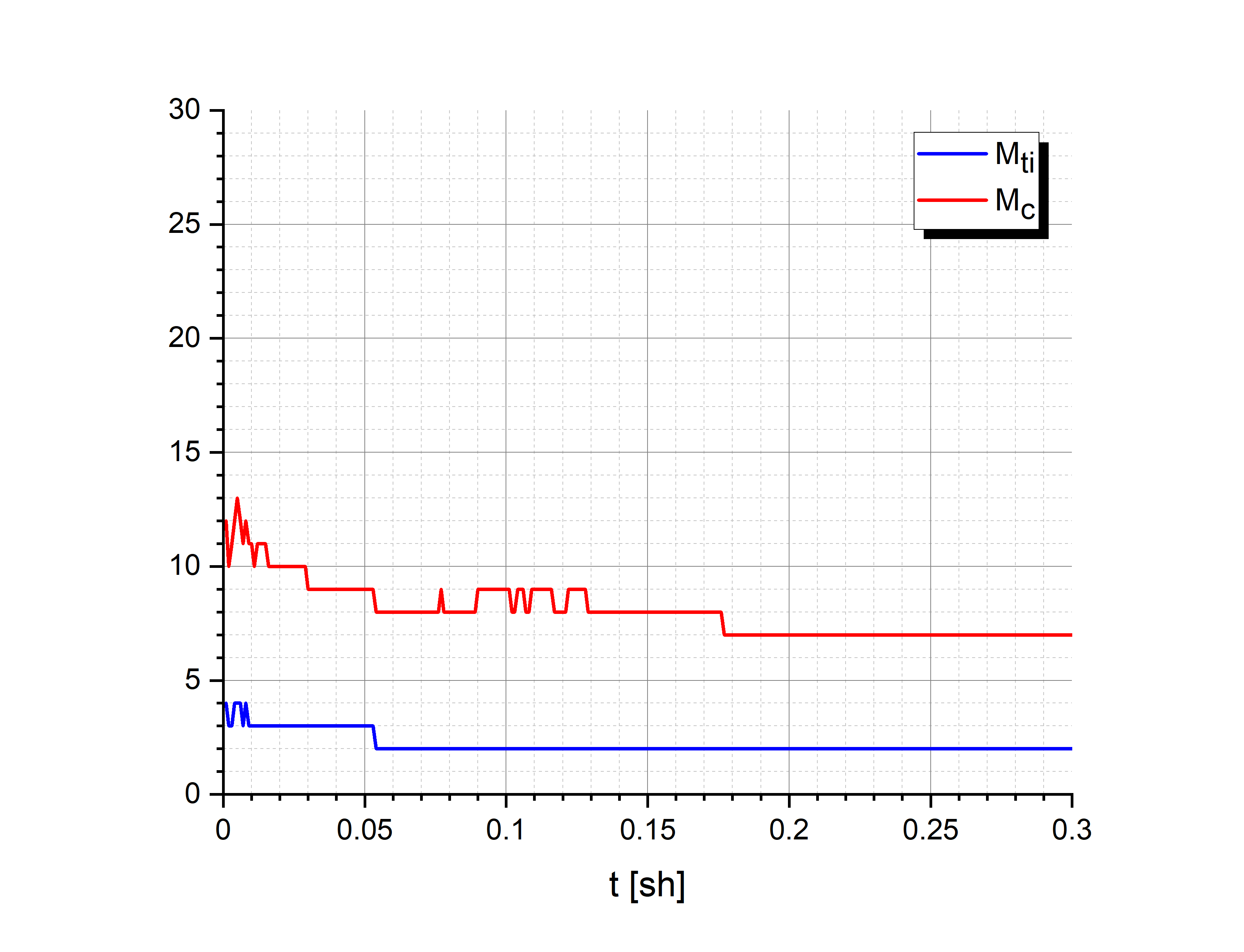}}
\hfill
\subfloat[$\ell_{max}$=4  \label{V4-1e-3}]
{\includegraphics[scale=0.28]{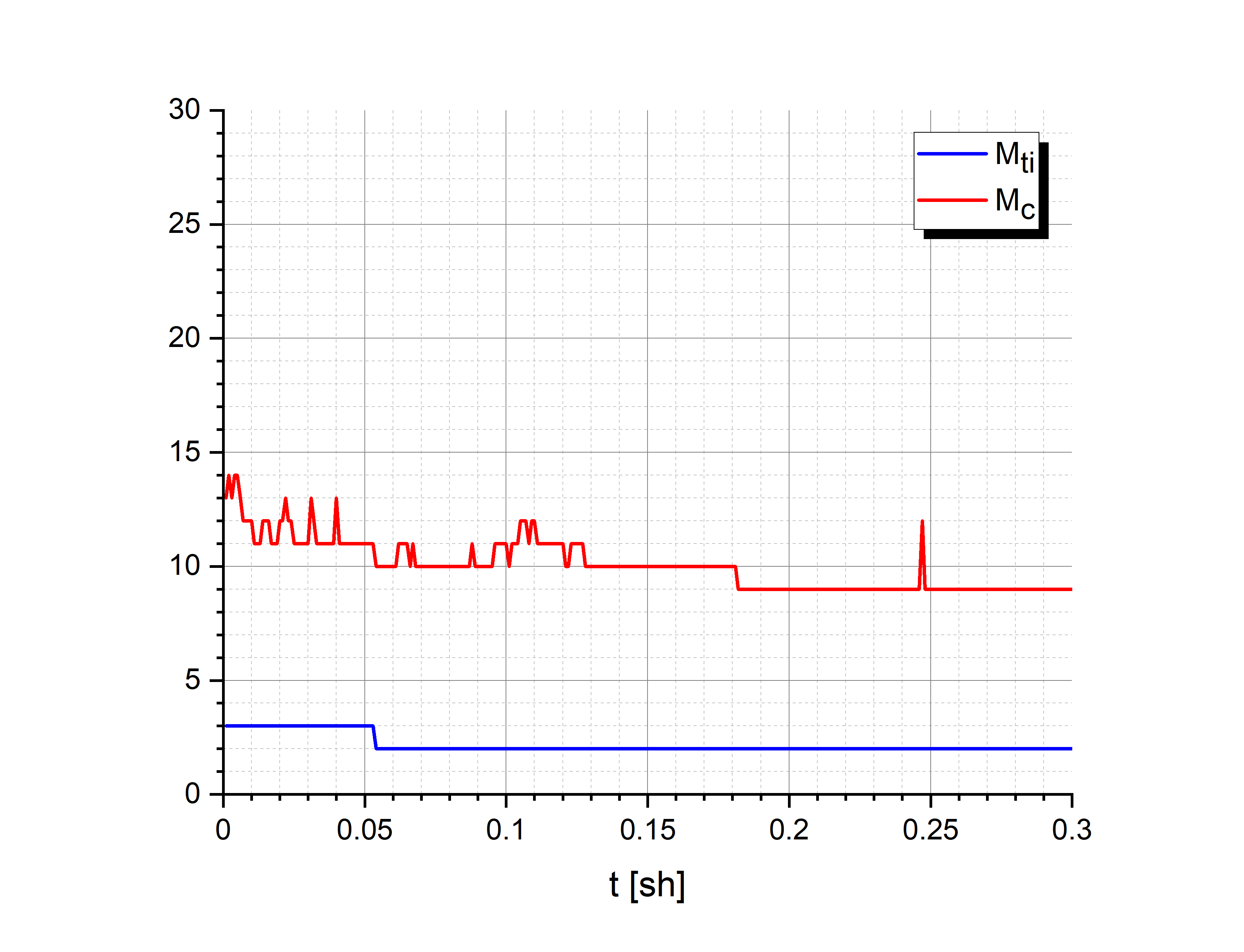}}
 \caption{\label{algorithm-256-001}  Number of transport iterations ($M_{ti}$) and number of cycles ($M_c$)   at each time step  in the FC test with 256 groups  and $\Delta t = 1 \! \times \! 10^{-3}$ sh  over  $t\in$ [0,  0.3\,sh].}
\end{figure}

\begin{figure}[h!]
\centering
\subfloat[$\ell_{max}$=3  \label{V3-2e-3}]
{\includegraphics[scale=0.28]{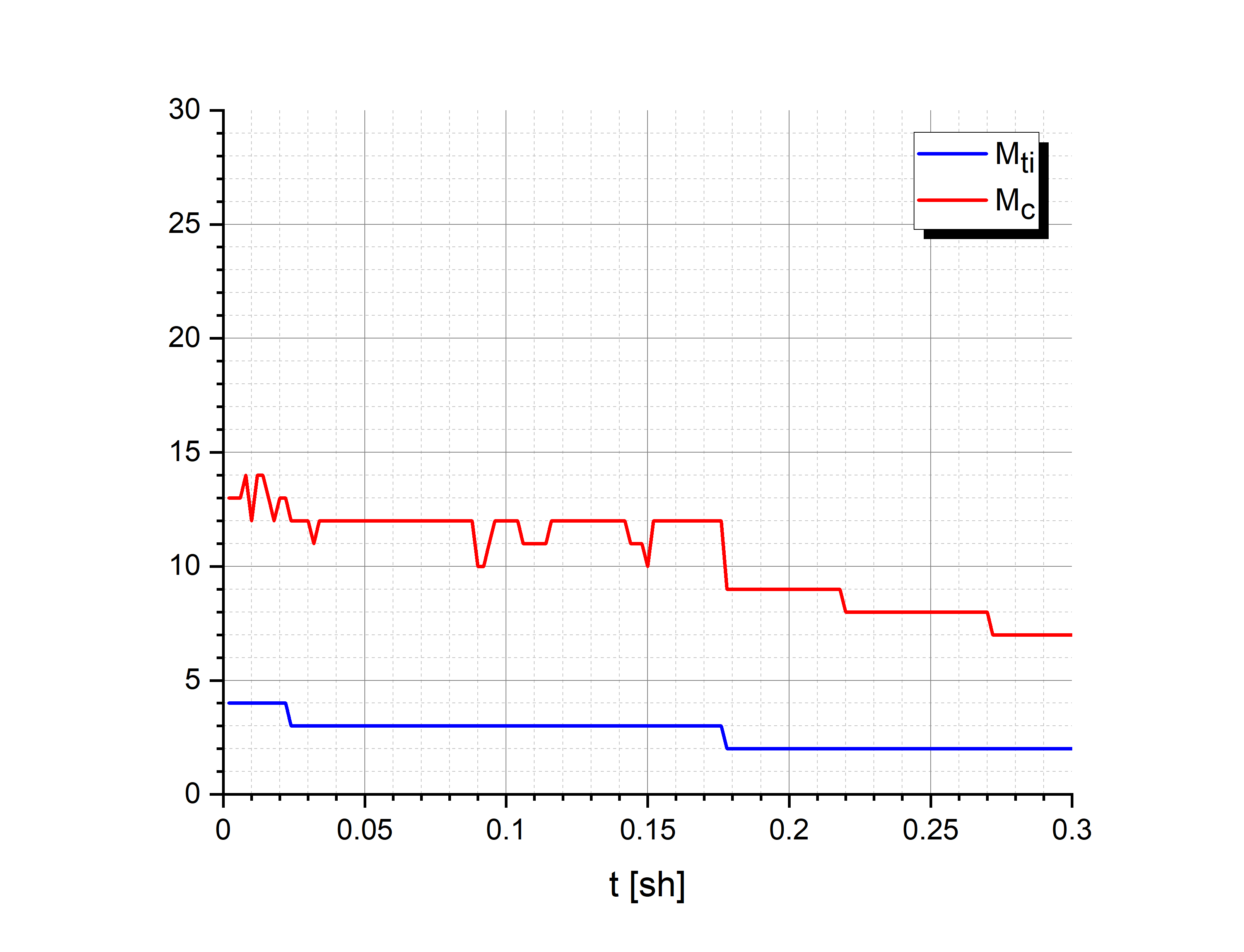}}
\hfill
\subfloat[$\ell_{max}$=4  \label{V4-2e-3}]
{\includegraphics[scale=0.28]{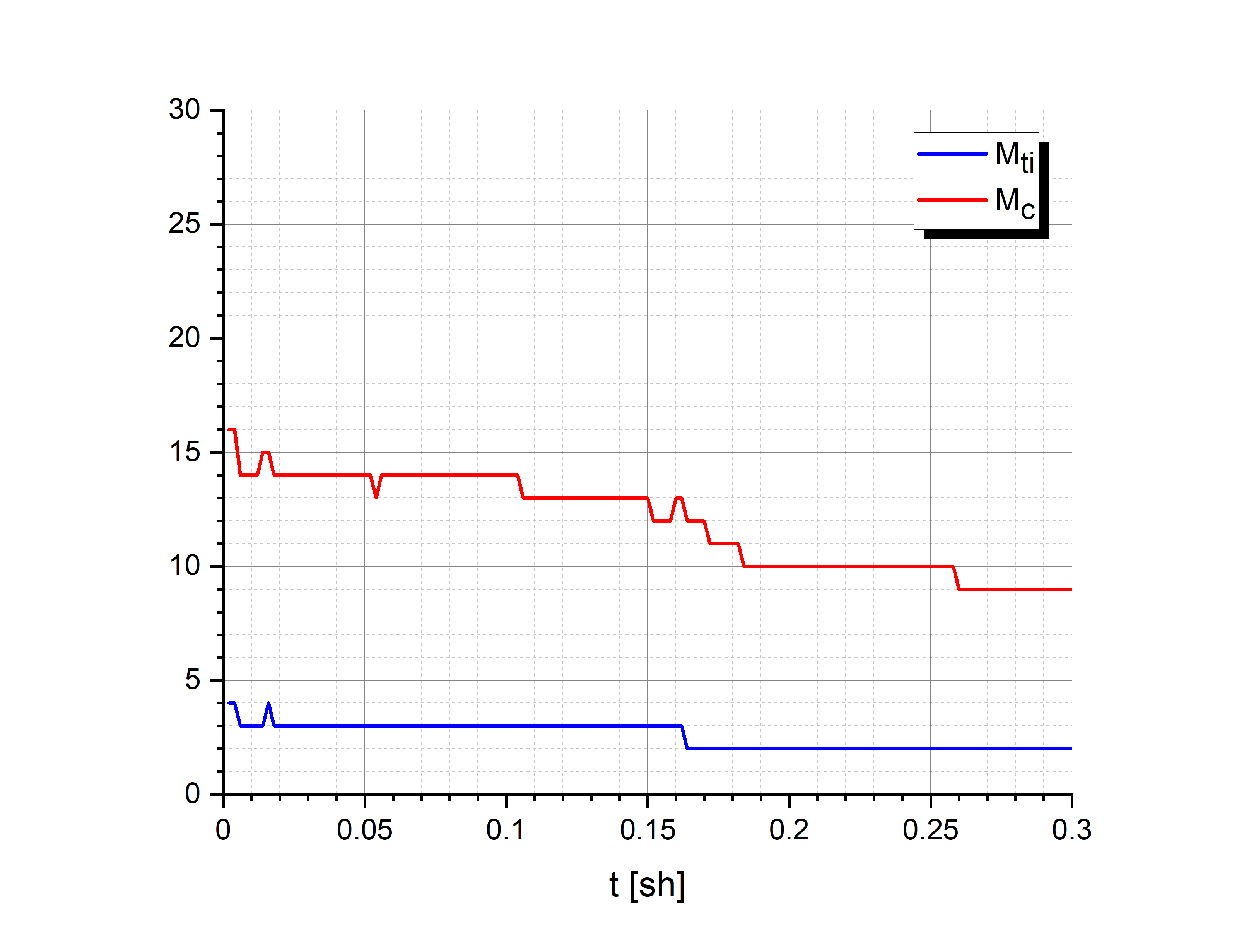}}
  \caption{\label{algorithm-256-002}  Number of transport iterations ($M_{ti}$) and number of cycles ($M_c$)   at each time step  in the FC test with 256 groups  and $\Delta t = 2 \! \times \! 10^{-3}$ sh  over  $t\in$ [0, 0.3\,sh].}
\end{figure}

\begin{figure}[h!]
\centering
\subfloat[$\ell_{max}$=4  \label{V4-4e-3}]
{\includegraphics[scale=0.28]{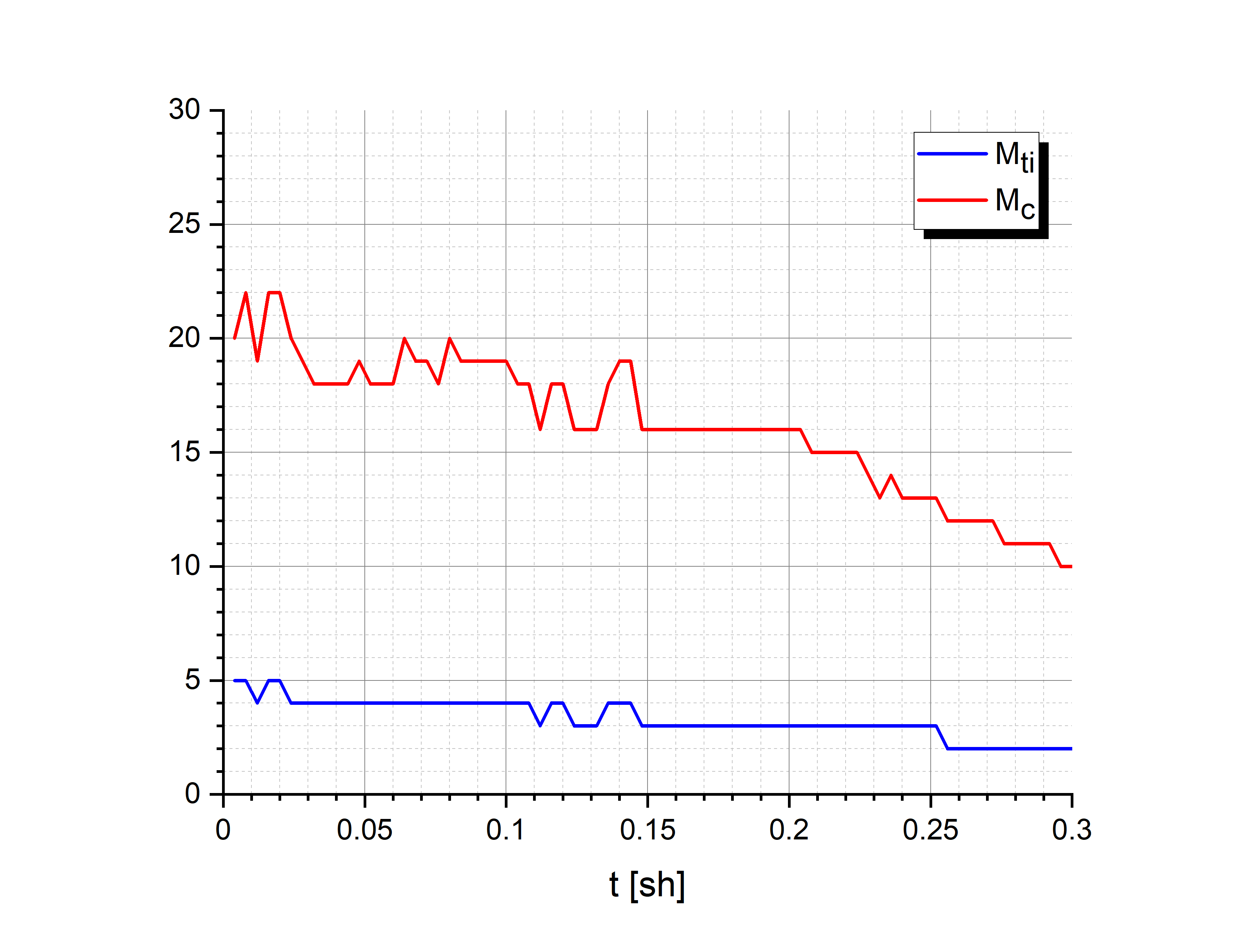}}
\hfill
\subfloat[$\ell_{max}$=5  \label{V5-4e-3}]
{\includegraphics[scale=0.28]{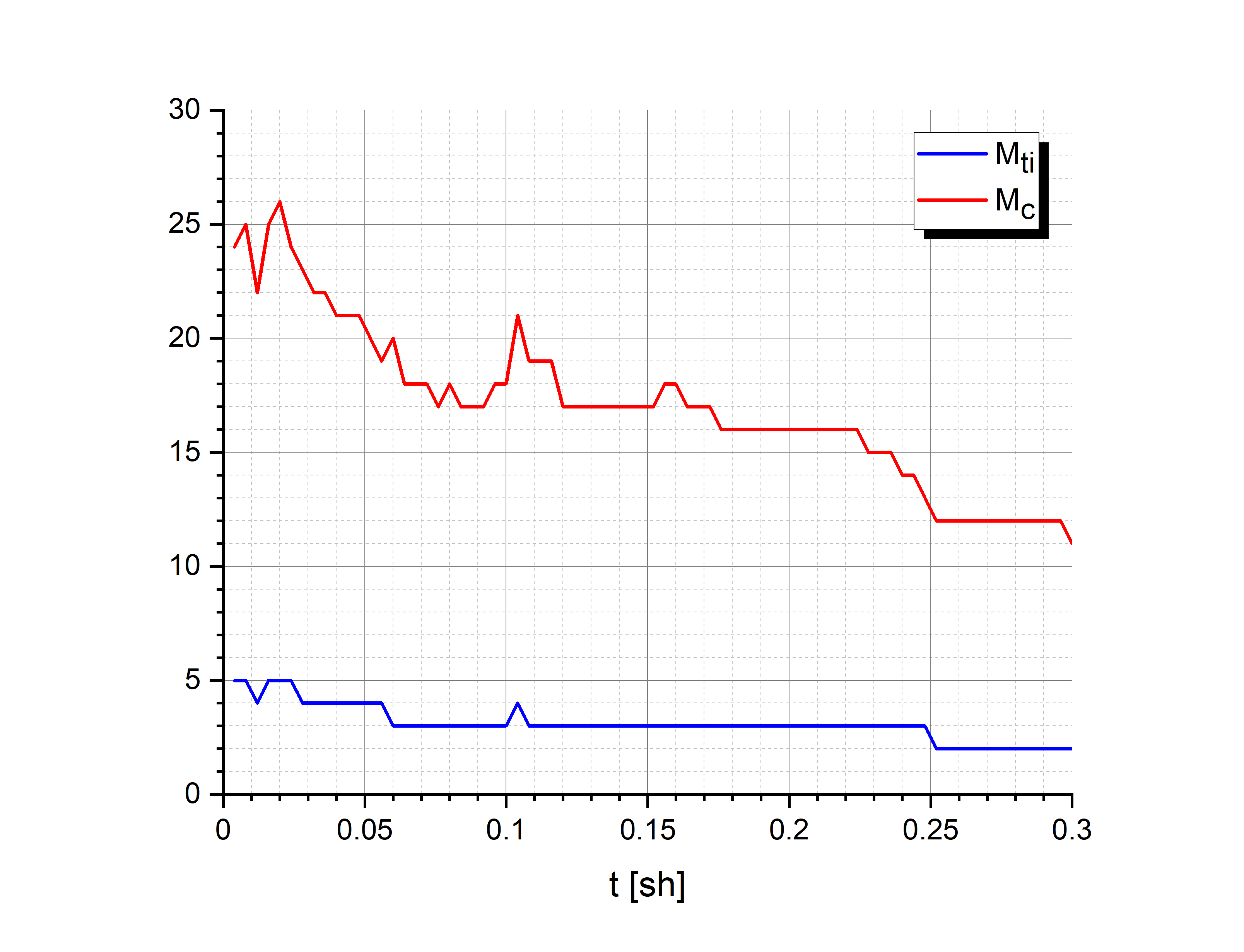}}
  \caption{\label{algorithm-256-004}  Number of transport iterations ($M_{ti}$) and number of cycles ($M_c$)   at each time step  in the FC test with 256 groups  and $\Delta t = 4 \! \times \! 10^{-3}$ sh  over  $t\in$ [0,  0.3\,sh].}
 \end{figure}

Figures  \ref{algorithm-256-001}-\ref{algorithm-256-004}  show extra details of
performance of algorithm with different  values of $\ell_{max}$.
The plots   present the total numbers of transport iterations, $M_{ti}$, and cycles, $M_c$, 
 at each  time step.
There are several different stages in   Marshak wave evolution in this test problem.  
At the initial stage, the rate of change of the solution is large.
 This leads to more transport iterations and cycles compared to
the next stage during which we observe a well developed wave. 
At the last stage, the  solution approaches a steady-state regime.
The number of transport iterations is very small.
Figures  \ref{algorithm-256-001}-\ref{algorithm-256-004} demonstrate these effects as well as  
influence of the time step value.

\begin{table}[h!]
	\centering
	\caption{\label{tbl-256-004} Performance of    algorithms  in the FC test with 256 groups using different time steps $\Delta t$}
\begin{tabular}{|c||c|c||c|c||c|c|}
  \hline
   &  \multicolumn{2}{c||}{    $  \Delta t  = 1 \! \times \! 10^{-2} \, sh $   }  & 
    \multicolumn{2}{c||}{  $ \Delta t= 2 \! \times \! 10^{-2} \, sh $  }   & \multicolumn{2}{c|}{  $ \Delta t= 4\! \times \! 10^{-2} \,  sh $ }  \\ \hline  
 $ \ell_{max} $          &   $N_{ti}$    &   $N_{c}$    & $N_{ti}$         &   $N_{c}$  & $ N_{ti}$     &   $N_{c}$        \\ \hline
2                & 1208          &  2970          & 653           &  1605   & 432           &  1014  \\ \hline
3                 & 658           &  2408          & 399           &  1554   & 300           &  1109    \\ \hline
4                 & 653           &  3010          & 384           &  1789   & 249           &  1215   \\ \hline 
5               & 653           &  3297          & 384           &  2023   & 232           &  1278   \\ \hline  
6                & 653           &  3509          & 384           &  2283   & 232           &  1490   \\ \hline    
7              & 653           &  3534          & 384           &  2443   & 232           &  1700  \\ \hline   
  \end{tabular}
\end{table}


\section{Conclusions \label{conclusion}}

We developed a new multilevel iteration method for solving multigroup TRT problems in 1D slab geometry.
The obtained results are promising. The algorithm is rapidly convergent.
 It is stable with independent discretization of the high-order RTE and low-order equations at each level.
 Efficiency of the algorithm can be optimized by choosing the number of $\textsf{V}$-cycles per outer (transport) iteration.
Further detailed analysis of iteration algorithm is necessary.
Future work will include development of  methods with the mixed type low-order
problems using multiple grids in photon energy and various multigrid cycles to further optimize performance of iteration schemes
and computational cost.  Various discretization schemes will be considered.
The proposed  multilevel  method can be extended to multidimensional geometries, and is suitable for computations on high-performance computing systems with GPUs.

\section*{Acknowledgements}

This manuscript has been authored by Lawrence Livermore National Security, LLC under Contract No. DE-AC52-07NA27344 with the US. Department of Energy. The United States Government retains, and the publisher, by accepting the article for publication, acknowledges that the United States Government retains a non-exclusive, paid-up, irrevocable, world-wide license to publish or reproduce the published form of this manuscript, or allow others to do so, for United States Government purposes. LLNL-CONF-2001051.

\bibliographystyle{elsarticle-num}
\bibliography{dya-tsh-arxiv-2024}

\end{document}